\documentclass[11pt,leqno]{article}
\usepackage{amssymb, amscd, amsmath, amsthm}

\allowdisplaybreaks

\def\ve{\varepsilon}

\def\mod{\,\text{\rm mod}\;}
\def\beq{\begin{equation}}
\def\eeq{\end{equation}}
\def\cite#1{{\rm [#1]}}

\newtheorem{theorem}{Theorem}

\theoremstyle{definition}

\newtheorem*{rema}{Remark}

\begin{document}

\numberwithin{equation}{section}
\title{On the singular series in the prime $k$-tuple conjecture}

\author{J\'anos Pintz\thanks{Supported by OTKA Grants K72731, K67676 and ERC-AdG.228005.}
}

\date{}
\maketitle
\setcounter{section}{1}

1. Gallagher's result \cite{Gal} about the mean-value of the singular series
\beq
\mathfrak S(\mathcal H) = \prod_p \left(1 - \frac{\nu_p}{p}\right) \left(1 - \frac1p\right)^{-k}
\label{eq:1.1}
\eeq
over all $k$-element sets
\beq
\mathcal H = \{h_i\}^k_{i = 1}, \quad h_i \text{ different}, \quad h_i \in [1, H]
\label{eq:1.2}
\eeq
(where $\nu_p = \nu_p(\mathcal H)$ denotes the number of residue classes occupied by $\mathcal H \mod p$)
played a crucial role in the proof of
\beq
\liminf_{n \to \infty} \frac{p_{n + 1} - p_n}{\log p_n} = 0
\qquad \text{ \cite{GPY}},
\label{eq:1.3}
\eeq
where $p_n$ denotes the $n$\textsuperscript{th} prime.

However, a brief study of the works \cite{GPY} or \cite{GMPY} shows that actually we need a far weaker result:
to show the existence of at least one admissible set $|\mathcal H| = k_\nu$ (i.e.\ $\nu_p < p$ for any prime $p$) for a series of $k_\nu \to \infty$ and a $C > 0$ absolute constant with the property
\beq
S_{\mathcal H}(H) := \frac1H \sum^H_{h = 1} \frac{\mathfrak S(\mathcal H \cup \{h\})}{\mathfrak S(\mathcal H)} \geq C \quad \text{ for }\ H > C_0 (|\mathcal H|).
\label{eq:1.4}
\eeq

In \cite{Pin} we showed in a very simple way that \eqref{eq:1.4} holds for a particular $\mathcal H_k$ for any $k$, even for all individual even~$h$ (for odd $h$ we had $\mathfrak S(\mathcal H \cup \{h\}) = 0$).
In the present note we prove beyond \eqref{eq:1.4} that for any $\mathcal H \subset [1, H]$ we have $S_{\mathcal H}(H) \to 1$ if $H \to \infty$, and investigate the rate of convergence depending on~$k$.
By induction on~$k$, this implies Gallagher's result \cite{Gal}
\beq
\sum_{|\mathcal H| = k,\ \mathcal H \subset [1,H]} \mathfrak S(\mathcal H) \sim H^k \quad \text{ for }\ H \to \infty,
\label{eq:1.5}
\eeq
however, $S_{\mathcal H}(H) \to 1$ yields more information about the singular series than the global average~\eqref{eq:1.5}.

\begin{theorem}
\label{th:1}
With the notation \eqref{eq:1.1}, \eqref{eq:1.2}, \eqref{eq:1.4} we have for any sufficiently small $\ve > 0$
\beq
S_{\mathcal H}(H) = 1 + O(\ve) \quad \text{ if } \ H \geq \exp (k^{1/\ve}).
\label{eq:1.6}
\eeq
\end{theorem}

\begin{rema}
The following relations can be proved in an even more simple way for every~$\mathcal H$.
\end{rema}

\noindent
{\bf Theorem 1'.}
{\it With the notation \eqref{eq:1.1}--\eqref{eq:1.2} we have
\beq
S_{\mathcal H}(H) \geq 1 + O(\ve)\ \ \text{ if }\ H \geq \exp \left(\frac1{\ve} \frac{k}{\log k}\right)
\label{eq:1.7}
\eeq
and with some absolute constants $c_1$, $c_2$}
\beq
S_{\mathcal H}(H) \geq c_1 \ \ \text{ {\it if} } \ H \geq \exp \left(c_2 \frac{k}{\log k}\right).
\label{eq:1.8}
\eeq

\begin{rema}
As mentioned in the introduction, for \eqref{eq:1.3}, but actually in most other applications for problems involving small gaps between primes and almost primes we need just lower estimates for the singular series.
This fact gives an additional significance for simple proofs of lower estimates for expressions involving the singular series.
\end{rema}

\bigskip
\setcounter{section}{2}
\setcounter{equation}{0}
2. Let us study first the ratio $\mathfrak S(\mathcal H')/\mathfrak S(\mathcal H)$ for a single $\mathcal H' := \mathcal H \cup \{h\}$, $h \notin \mathcal H$, with the notation
\beq
\nu'_p := \nu_p(\mathcal H'), \quad
y := \frac{5\log H}{6}, \quad
P := \prod_{p \leq y} p, \quad
\Delta := \prod^k_{i = 1} (h - h_i).
\label{eq:2.1}
\eeq
Then we have
\beq
\frac{\mathfrak S(\mathcal H')}{\mathfrak S(\mathcal H)} = \prod_p \frac{1 - \frac{\nu'_p}{p}}{\left(1 - \frac{\nu_p}{p}\right) \left(1 - \frac1p\right)} := \underset{p\leq y}{\prod\nolimits_1} \cdot
\underset{{\scriptstyle p > y\atop \scriptstyle p \mid \Delta}}{\prod\nolimits_2} \cdot
\underset{{\scriptstyle p > y\atop \scriptstyle p \nmid \Delta}}{\prod\nolimits_3} .
\label{eq:2.2}
\eeq
For $p\nmid \Delta$ we have $\nu'_p = \nu_p + 1$, otherwise $\nu'_p = \nu_p$, hence
\begin{align}
\prod\nolimits_3 &= \prod_{p > y} \left(1 + O\left(\frac{k}{p^2}\right)\right) = 1 + O \left( \frac{k}{y \log y} \right),
\label{eq:2.3}\\
\prod\nolimits_2 &= \prod_{p > y,\ p \mid  \Delta} \left(1 - \frac1p\right)^{-1} = \exp \left(O\left(\sum_{p > y, \ p \mid \Delta} \frac1{p}\right)\right) .
\label{eq:2.4}
\end{align}
Since by the Prime Number Theorem
\[
\sum\limits_{p \mid \Delta} \log p \leq  \log \Delta \leq 2 ky,
\]
the sum over $1/p$ in the error term is maximal if the relevant primes with $p \mid \Delta$ are as small as possible satisfying the inequality $\sum\limits_{p \in \mathcal P^*} \log p \leq 2ky$, consequently
\beq
\sum_{p\mid \Delta, \ p > y} \frac1p \leq \sum_{y < p < 4ky} \frac1p \leq \log \frac{\log(5ky)}{\log y} \leq 2\ve.
\label{eq:2.5}
\eeq

\eqref{eq:2.3}--\eqref{eq:2.5} imply that we have for any individual $h$
\beq
\prod\nolimits_2 \prod\nolimits_3 (h) = 1 + O(\ve).
\label{eq:2.6}
\eeq
We will now study the average of $\prod_1(h)$ for $h \leq H = M P + r$, $0 \leq r < P$, $M \to \infty$.
Due to the periodicity of $\prod_1(h)$ (with period $P$) it is sufficient to average over $h \in [1,P]$.
For any $p \leq y$ we have exactly $\nu_p$ possibilities $\mod p$ for $h$ with $\nu'_p = \nu_p$ and $p - \nu_p$ with $\nu'_p = \nu_p + 1$.
Consequently,
\beq
\frac1{P} \sum^P_{h = 1} \prod\nolimits_1(h) = \prod_{p \mid P} \frac{\left\{\frac{\nu_p}{p} \left(1 - \frac{\nu_p}{p}\right) + \left(1 - \frac{\nu_p}{p}\right) \left(1 - \frac{\nu_p + 1}{p} \right)\right\}}{\left(1 - \frac{\nu_p}{p}\right) \left( 1 - \frac1{p}\right)} = 1.
\label{eq:2.7}
\eeq

Formulas \eqref{eq:2.2} and \eqref{eq:2.6}--\eqref{eq:2.7} prove \eqref{eq:1.6}.

\medskip
To prove Theorem~1' we can replace \eqref{eq:2.4}--\eqref{eq:2.5} with the trivial relation $\Pi_2 \geq 1$, thereby obtaining
\beq
\Pi_2 \Pi_3(h) \geq \begin{cases}
1 + O(\ve) &\text{if }\ \frac{k}{\log H \log_2 H} = O(\ve),\\
c_1 &\text{if }\ \frac{k}{\log H \log_2 H} \leq c_3.
\end{cases}
\label{eq:2.8}
\eeq
This, together with \eqref{eq:2.7} proves Theorem~1'.

\bigskip

\noindent
{\small J\'anos {\sc Pintz}\\
R\'enyi Mathematical Institute of the Hungarian Academy
of Sciences\\
Budapest\\
Re\'altanoda u. 13--15\\
H-1053 Hungary\\
E-mail: pintz@renyi.hu}

\end{document}